\documentclass[leqno,11pt]{amsart}
\usepackage{fullpage}
\usepackage[english]{babel}
\usepackage{amsmath,amsthm,amssymb,amsfonts,mathtools,mathrsfs}
\usepackage[alphabetic]{amsrefs}
\usepackage[frak=pxtx,cal=euler,scr=boondoxupr]{mathalfa}	

\usepackage{xcolor}	

\newcommand{\spacesymbol}[1]{\mathbb{#1}}
\newcommand{\sans}[1]{\mathsf{#1}}
\newcommand{\R}{\spacesymbol{R}}
\newcommand{\C}{\spacesymbol{C}}
\newcommand{\N}{\spacesymbol{N}}

\newcommand{\m}[2]{\left\langle #1, #2 \right\rangle}
\DeclareMathOperator\tr{tr}

\DeclareMathOperator\dive{div}
\DeclareMathOperator\Ric{Ric}
\DeclareMathOperator\ric{ric}

\DeclareMathOperator\Aut{Aut}
\DeclareMathOperator\Der{Der}

\newtheorem{maintheorem}{Theorem}

\newtheorem{theorem}{Theorem}[section]

\newtheorem{prop}[theorem]{Proposition}
\newtheorem{corollary}[theorem]{Corollary}
\newtheorem{lemma}[theorem]{Lemma}

\theoremstyle{definition}
\newtheorem{remark}[theorem]{Remark}
\newtheorem*{notation}{Notation}
\newtheorem*{acknowledgements}{Acknowledgements}

\theoremstyle{remark}
\newtheorem{example}[theorem]{Example}

\title{Inhomogeneous deformations of Einstein Solvmanifolds}
\author[A. Thompson]{Adam Thompson}
\address{School of Mathematics and Physics, The University of Queensland, St Lucia QLD 4072, Australia}
\email{adam.thompson@uq.net.au}

\begin{document}

\maketitle

\begin{abstract}
For each non-flat, unimodular Ricci soliton solvmanifold \((\sans{S}_0,g_0)\), we construct a one-parameter family of complete, expanding, gradient Ricci solitons that admit a cohomogeneity one isometric action by \(\sans{S}_0\). The orbits of this action are hypersurfaces homothetic to \((\sans{S}_0,g_0)\). These metrics are asymptotic at one end to an Einstein solvmanifold. In the one-parameter family, exactly one metric is Einstein, and exactly one has orbits that are isometric to \((\sans{S}_0,g_0)\).
\end{abstract}

\section{Introduction}

A Ricci soliton is a Riemannian manifold \((M,g)\) whose Ricci curvature satisfies \begin{equation}\label{Eq:RS}
\ric(g) = \lambda g+ \frac{1}{2}\mathcal{L}_Xg,\quad \lambda\in \R,
\end{equation} where \(X\in \mathfrak{X}(M)\) is a complete vector field. In this article, we are interested in expanding solitons: \(\lambda<0\). If the vector field \(X\) is a Killing field (i.e. \(\mathcal{L}_Xg=0\)), then \(\ric(g)=\lambda g\), so \(g\) is an \emph{Einstein metric}. If \(X=\nabla f\) for some \(f\in C^\infty(M)\), we say the soliton is of \emph{gradient type}.

A driving force in the study of Einstein metrics and Ricci solitons with symmetry has been the Alekseevsky conjecture: that a homogeneous Einstein manifold of negative scalar curvature is diffeomorphic to \(\R^n\). This has lead to many structural results for Einstein metrics and Ricci solitons with non-compact symmetry, even beyond the homogeneous setting, cumulating in its recent resolution \cite{BohmLafuente23}. Despite these structural results, there are very few examples of complete Einstein manifolds and Ricci solitons with non-compact symmetry beyond the homogeneous setting. As far as we know, all known examples are: the expanding solitons on vector bundles constructed in \cite{Lott07}, which have arbitrary codimension, and are characterised in terms of a harmonic map equation and an equation on the base manifold; the \((S^1\times \sans{N})\)-invariant Einstein metrics in \cite{Hervik08}, where \(\sans{N}\) is a Ricci soliton nilmanifold; the three-dimensional expanding soliton from \cite{Ramos13} lifts to an \(\R^2\)-invariant soliton on \(\R^3\); the quaternionic K{\"ahler} metrics in \cite{CortesSahaThung21,CortesSaha22}; and the expanding K{\"ahler}-Ricci solitons in \cite{MaschlerReam22} which are invariant under a Heisenberg group action. Note that all known examples of non-trivial Ricci solitons with non-compact symmetry have either an abelian symmetry group or special holonomy.

In this article, we construct infinite families of complete Ricci solitons that are invariant under solvable Lie groups whose orbits are hypersurfaces. These generalise most of the known examples, and most of them do not have abelian symmetry or special holonomy.

\begin{maintheorem}\label{Thm:MainThm}
Let \((\sans{S}_0,g_0)\) be a unimodular solvsoliton with negative scalar curvature.  There exists a one-parameter family of pairwise non-homothetic \(\sans{S}_0\)-invariant, complete, gradient, expanding Ricci soliton metrics on \(M=\R\times \sans{S}_0\) that are not homogeneous. The \(\sans{S}_0\)-orbits of each of these Ricci solitons are hypersurfaces homothetic to \((\sans{S}_0,g_0)\). Exactly one of the metrics in the family is Einstein.
\end{maintheorem}

A solvsoliton is a simply connected solvable Lie group with a left-invariant metric, \(g_0\), whose Ricci operator satisfies\[\Ric(g_0)=\lambda_0I+D,\qquad \lambda_0<0,\quad D\in \Der(\mathfrak{s}_0),\] where \(\mathfrak{s}_0=\operatorname{Lie}(\sans{S}_0)\). Note that solvsolitons are Ricci solitons \cite[Section~2]{Lauret11}.

Left invariant metrics on solvable Lie groups always have non-positive scalar curvature. If \(\sans{S}_0\) is scalar flat, then the construction in Theorem~\ref{Thm:MainThm} is still valid, but no longer yields a one-parameter family. Instead, we obtain a single soliton. We note however that, in this case, our Ansatz reduces to a warped product Ansatz and the existence of the soliton follows from the analysis in \cite{Ramos13}. See \S \ref{Section:NoScal} for a discussion of this case.

Given a soliton from Theorem~\ref{Thm:MainThm}, the soliton vector field, \(X=\nabla f\), is \(\sans{S}_0\)-invariant. Therefore, if \(\Gamma\le \sans{S}_0\) is a cocompact lattice, the quotient \(\Gamma\backslash M=\R\times (\Gamma\backslash \sans{S}_0)\) is a gradient soliton that has two non-compact ends.

\begin{example}
Let \(N^3=\Gamma\backslash\mathrm{Nil}\) be a three-manifold with Nil-geometry, where \(\mathrm{Nil}\) is the Heisenberg group. Previously there were two known (inhomogeneous) complete soliton metrics on \(\R\times N\) invariant under the \(\mathrm{Nil}\) action: the K\"ahler-Ricci soliton in \cite{MaschlerReam22}, and the quaternionic K{\"a}hler metric in \cite{CortesSaha22}. Theorem~\ref{Thm:MainThm} gives a one-parameter family of metrics which includes these two known examples.
\end{example}

Note that the higher-dimensional examples constructed in \cite{MaschlerReam22} also appear in Theorem~\ref{Thm:MainThm} if we take \(\sans{S}_0\) to be the \((2n+1)\)-dimensional generalised Heisenberg group.

The solitons in Theorem~\ref{Thm:MainThm} are \emph{cohomogeneity one}: the \(\sans{S}_0\)-orbits have dimension one less than \(M\). Many examples in Riemannian geometry have been found using a cohomogeneity one Ansatz; for example, Einstein metrics on low-dimensional spheres constructed by B{\"o}hm \cite{Boehm98}, and the first complete inhomogeneous examples of nearly K{\"ahler} 6-manifolds constructed by Foscolo-Haskins \cite{FoscoloHaskins17}. Also, many of the early examples of Ricci solitons, such as the Cigar soliton and Bryant solitons, are invariant under a cohomogeneity one action by a compact group; see \cite[Chapter~1]{Chow07}. Many examples of gradient Ricci solitons that admit a cohomogeneity one action by a compact group have also been constructed over the past 15 years; see \cite{DancerWang10,BuzanoDancerWang15I,BuzanoDancerWang15II,Appleton18,Buttsworth21,Wink20,Wink21}, together with the references therein.

On the other hand, there are far fewer examples in the case that the group \(\sans{G}\) acting by cohomogeneity one is non-compact, as we have mentioned above.

If \(\sans{S}\) is a solvable Lie group equipped with a left-invariant Einstein metric, \(g^S\), then \(\sans{S}\) admits a codimension 1, unimodular, closed subgroup \(\sans{S}_0\) which is a solvsoliton with the induced metric \cite{LafuenteLauret14}. Conversely, if \((\sans{S}_0,g_0)\) is a solvsoliton, where \(\sans{S}_0\) is a unimodular, solvable Lie group, then \(M=\R\times \sans{S}_0\) admits a homogeneous Einstein metric extending \(g_0\) \cite[Proposition~6.1]{LafuenteLauret14}. The metrics in Theorem~\ref{Thm:MainThm} are obtained by extending \(\sans{S}_0\) in an inhomogeneous way. In \cite{AlekseevskyNikolayevsky21}, a different approach to generalising the connection between \(\sans{S}\) and \(\sans{S}_0\) to an inhomogeneous setting is considered.

The previously known examples that appear in our construction also admit quotients that have two non-compact ends. These are asymptotic at one end to (a quotient of) an Einstein solvmanifold. Moreover, the Einstein manifold in \cite{CortesSaha22} is asymptotic at its other end to \(\R H^4\); the K\"ahler-Ricci solitons in \cite{MaschlerReam22} are asymptotic to a cone whose link is \(\Gamma\backslash \sans{S}_0\) at their other end; and, the solitons in \cite{Ramos13} as asymptotic to a cone whose link is the two-torus.  

In is natural to expect that the family we produce has similar asymptotics to these. This is indeed the case.

\begin{maintheorem}\label{Thm:Asymptotics}
Let \((\sans{S}_0,g_0)\) be a unimodular solvsoliton with negative scalar curvature, \(\Gamma\le \sans{S}_0\) a cocompact lattice, and let \((M,g)\) be a Ricci soliton from Theorem~\ref{Thm:MainThm}. Then, \((\Gamma\backslash M,g)\) has two non-compact ends. At one end, \((\Gamma\backslash M,g)\) is asymptotic to (a quotient of) the Einstein solvmanifold obtained from \((\sans{S}_0,g_0)\) by a one-dimensional extension; at the other end, \((\Gamma\backslash M,g)\) has one of the three following asymptotics:

\begin{enumerate}
\item If \((\Gamma\backslash M,g)\) is Einstein, then it is asymptotic to a space of constant negative curvature;

\item If the \(\sans{S}_0\)-orbits in \((\Gamma\backslash M,g)\) are all isometric, then it is asymptotic at its other end to the metric \[ds^2 + e^{\frac{1}{-\lambda}(\log s) D_*}g_0,\] where \(\lambda\) is the cosmological constant of \((\Gamma\backslash M,g)\), and \(D\in \Der(\mathfrak{s}_0)\) is the soliton derivation;

\item If \((\Gamma\backslash M,g)\) is a non-trivial soliton (i.e. not Einstein) and the \(\sans{S}_0\)-orbits have different scalar curvature, then it is asymptotically conical.
\end{enumerate}
Moreover, if \((M,g)\) is a non-trivial soliton, then the potential function \(f:M\to\R\) is convex (i.e. has positive definite Hessian).
\end{maintheorem}

Let us outline the construction used in Theorem~\ref{Thm:MainThm}. An \(\sans{S}_0\)-invariant metric on \(M=\R\times \sans{S}_0\) can be described by a one-parameter family of inner products \(\{g(s)\}\) on \(\mathfrak{s}_0=\operatorname{Lie}(\sans{S}_0)\). We make the following Ansatz: for all \(s\in \R\), the inner product \(g(s)\) can be written as \[g(s) = c^2(s) e^{h(s)D*}g_0 = g_0\big( c^2(s)e^{2h(s)D}\cdot,\cdot\big),\] where \(c:\R\to \R_+\), \(h:\R\to \R\) are smooth functions, and \(D\in \Der(\mathfrak{s}_0)\) is the soliton derivation (that is, \(\Ric(g_0)=\lambda_0I+D\) for some \(\lambda_0\in \R\)). 

The equation \eqref{Eq:RS} becomes an ordinary differential equation (assuming the soliton potential is also \(\sans{S}_0\)-invariant) for the shape operator \(L(s)= \frac{1}{2}g'(s)g(s)^{-1}\) of the \(\sans{S}_0\)-orbits and for the derivative of the soliton potential \(f\). (We use a prime to denote differentiation with respect to \(s\), \(\cdot'=d/ds\).) Existence and uniqueness of ODEs immediately gives us short time solutions to the equations defining the cohomogeneity one solitons; the non-trivial task is to prove that we have solutions that correspond to complete metrics. Since we are interested in the case that the orbit space is \(\R\), this will be true if and only if a solution is defined for all \(s\in \R\).  

Any unimodular solvsoliton \((\sans{S}_0,g_0)\) admits a one-dimensional Einstein extension \cite{LafuenteLauret14}. The Einstein manifold \((M,g^S)\) associated to our background solvsoliton \((\sans{S}_0,g_0)\) manifests in the ODE system as a stationary solution, which we denote \(\gamma^S\). Linearising the ODE at \(\gamma^S\) shows there is a two-dimensional unstable manifold of the system passing through \(\gamma^S\). Any integral curve beginning in the unstable manifold necessarily exists for all \(s\le 0\), and must converge to the stationary point \(\gamma^S\). We show that for a carefully chosen initial condition in the unstable manifold, the integral curve exists for all forwards times as well.

The paper has the following structure. In \S\ref{Section:Preliminaries} we explain our set up and recall the equations that define a cohomogeneity one Ricci soliton. \S\ref{Section:GeneralCase} is devoted to establishing the existence of the gradient Ricci solitons and computing their asymptotics. We then consider the Einstein case, and the case where the orbits are all isometric, in \S\ref{Section:EinsteinCase} and \S\ref{Section:NoScal} respectively. In Appendix~A we give some calculations that are omitted from \S\ref{Section:GeneralCase} and in Appendix~\ref{Appendix:ProperGroupActions} we prove some facts about proper group actions that are relevant to the article.

\begin{acknowledgements}
I would like to my supervisor, Ramiro Lafuente, for his support and advice. I would also like to thank Tim Buttsworth for his careful reading of a draft of this paper, and for his many helpful comments. The author was supported by an Australian Research Training Program scholarship.
\end{acknowledgements}

\section{Set up of the Cohomogeneity One Equations}\label{Section:Preliminaries}

Let \(\sans{S}_0\) be a simply connected, unimodular solvable Lie group that admits a solvsoliton metric, \(g_0\), with negative scalar curvature. Assume that \(\sans{S}_0\) acts properly, isometrically, and by cohomogeneity one on a Riemannian manifold \((M^{n+1},g)\). The action is necessarily free: the isotropy is a compact subgroup of \(\sans{S}_0\), and hence must be trivial.  Thus, the orbit space \(M/\sans{S}_0\) is either \(S^1\) or \(\R\). In fact, if \((M,g)\) is an expanding gradient Ricci soliton with an \(\sans{S}_0\)-invariant potential, then \(M/\sans{S}_0\) is necessarily \(\R\); this follows from \cite[Theorem~2]{PigaloRimoldiSetti11} and \cite[Theorem~F]{BohmLafuente23}. Since we are interested in expanding gradient Ricci solitons, we assume that \(M/\sans{S}_0=\R\).
 
 Choose a unit-speed geodesic \(\alpha:\R\to M^{n+1}\) meeting all of the orbits orthogonally; this gives an \(\sans{S}_0\)-equivariant diffeomorphism \(\R\times \sans{S}_0\simeq M\) by \((s,x)\mapsto x\cdot \alpha(s)\). Let \(N\) be the unit normal vector field of the orbits such that \(N_{\alpha(s)}=\alpha'(s)\) for all \(s\in \R\). The pullback of \(g\) to \(\R\times\sans{S}_0\) is given by \[g=ds^2 +g(s),\] where \(g(s)\) is a one-parameter family of left invariant metrics on \(\sans{S}_0\).  Note that while \(M\) is topologically a product, the geometry of \((M,g)\) is not a product. Using left invariance, and the background metric \(g_0\), we think of \(g(s)\) as a one-parameter family of endomorphisms of \(\mathfrak{s}_0=\operatorname{Lie}(\sans{S}_0)\): if \(\langle\cdot,\cdot\rangle_s\) denotes the inner product on \(\mathfrak{s}_0\) given by restricting \(g(s)\) to \(T_e\sans{S}_0\) then \(g(s):\mathfrak{s}_0\to \mathfrak{s}_0\) is defined by \[\langle \cdot,\cdot \rangle_s= \langle g(s)\cdot,\cdot\rangle_0.\] 

Let \(r(s):\mathfrak{s}_0\to \mathfrak{s}_0\) be the endomorphism of \(\mathfrak{s}_0\) associated to the Ricci tensor of \(\big(\sans{S}_0,g(s)\big)\): \[\ric (g(s))_e = g(s)_e\big(r(s)\cdot,\cdot).\] We also write \(r_0\) for the endomorphism associated to \(\ric(g_0)\). Observe that \(r_0=\lambda_0I+D\) for some \(\lambda_0<0\) and \(D\in \Der(\mathfrak{s}_0)\) since \((\sans{S}_0,g_0)\) is a solvsoliton. Note that the soliton derivation, \(D\), has non-negative eigenvalues; this is not true for generic solvsolitons, but it is true for unimodular solvsolitons (see the expression for \(D\) given in \cite[Proposition~4.3]{Lauret11}).

The Ricci operator, \(\Ric(g)\), of \((M,g)\) at \((s,x)\in \R\times \sans{S}_0\) can be expressed in terms of the shape operator \(L(s)\) of \(\{s\}\times \sans{S}_0\subset M\), and \(r(s)\); see for example \cite[Section~2]{EschenburgWang00}.  Equation \eqref{Eq:RS} with \(X=\nabla f=f' N\) then becomes the following ODE system for the shape operator \(L=\frac{1}{2}g'g^{-1}\) and \(f'\):

\begin{subequations}\label{Eq:Cohom1ODEs}
\begin{align}
\label{Eq:Cohom1ODEs(a)}&L'(s)+f'(s)L(s) = r(s) - \big(\tr L(s) \big)L(s)-\lambda I,\\
\label{Eq:Cohom1ODEs(b)}&\tr L'(s) +f''(s) = -\lambda - \tr \big(L^2(s)\big),\\
\label{Eq:Cohom1ODEs(c)}&{\ric}(g)\big(N,X)=0,\quad \text{for all }X \in T_p(\sans{S}_0\cdot p).
\end{align}\end{subequations}

We are interested in the Ansatz that the endomorphisms \(g(s):\mathfrak{s}_0\to \mathfrak{s}_0\) are of the form \begin{equation}\label{Eq:Ansatz}
g(s) = c(s)^2 e^{2h(s)D},
\end{equation} where \(c:\R\to \R_+\), \(h:\R\to \R\) are functions, and \(D\in \Der(\mathfrak{s}_0)\) is the soliton derivation; geometrically, this means that the orbits of \(\sans{S}_0\) are all homothetic to the fixed solvsoliton \((\sans{S}_0,g_0)\).

Recall that \(\Ric(c^2\varphi^*g_0)=c^{-2} \varphi^{-1}\Ric(g_0)\varphi\) for \(c>0\) and \(\varphi\in \Aut(\sans{S}_0)\).  Hence,  \[r(s)=  c(s)^{-2} e^{-h(s)D} \Ric(g_0)e^{h(s)D} = c(s)^{-2}r_0,\] since \(r_0\) and \(e^{tD}\) commute for all \(t\in \R\).

Moreover, we can calculate \[L(s)= \frac{1}{2}g^{-1}(s)g'(s) = \frac{1}{2}c(s)^{-2}e^{-2h(s)D} \bigg(2c'ce^{2hD} +2c^{2}h'e^{2hD}D\bigg) = \frac{c'}{c}I + h'D.\]  
Let us define functions \(x,y:\R\to \R\) by \[L(s) = x(s)I+y(s)D_0.\]
Here \(D_0\) is the trace-free part of \(D\). We can write \(x,y\) explicitly in terms of \(c,h\) as \[x=\m{L}{I}= \bigg( \frac{c'}{c} \bigg) + \frac{\tr D}{n}h', \quad y =\frac{\m{L}{D_0}}{\tr D_0^2}=h'. \] 
The variable \(x\) is (up to the factor \(n\)) the mean curvature of the \(\sans{S}_0\)-orbits (i.e. \(\tr L = nx\)); the variable \(y\) is measuring how much \(g(s)\) is moving in the \(\Aut(\mathfrak{s}_0)\)-direction.

We also let \(z(s)\) be the scalar curvature of the metric \(g(s)\) (equivalently, \(z\) is the scalar curvature of the orbits restricted to the orthogonal geodesic \(\alpha\)). Explicitly, \(z\) is given by \(z={s_0} c^{-2}\) where \(s_0\) is the scalar curvature of \((\sans{S}_0,g_0)\); note that \(s_0\) and \(\lambda_0\) are determined by the soliton derivation:\[\lambda_0 = \frac{-\tr D^2}{\tr D},\qquad s_0 = \frac{-\tr D^2}{\tr D}n+\tr D.\]

Finally, we denote by \(w\) the so-called dilation mean curvature \cite{BuzanoDancerWang15I} \[w=\tr L+f' =nx+f'.\]  Recall that the mean curvature is the logarithmic derivative of the ``volume'', \(\sqrt{\det g(s)}\), of the orbits. The dilation mean curvature is the logarithmic derivative of the ``dilated volume'' \(e^{f(s)}\sqrt{\det g(s)}\).
In these variables, the system \eqref{Eq:Cohom1ODEs} becomes
\begin{subequations}\label{Eq:RicciSolitonFlow}
\begin{align}
\label{Eq:RicciSolitonFlow(a)} x' &= \frac{z}{n}-\lambda-wx\\
\label{Eq:RicciSolitonFlow(b)} y'&=\frac{z}{s_0}-wy\\
\label{Eq:RicciSolitonFlow(c)} z'& = 2z\bigg( \frac{\tr D}{n}y-x\bigg)\\
\label{Eq:RicciSolitonFlow(d)} w'& =-\lambda -x^2n-y^2\tr D_0^2.
\end{align}
\end{subequations}

\begin{notation}
We denote points in \(\R^4\) by \(\gamma=(x,y,z,w)\).  We define the vector field \(F:\R^4\to \R^4\) to be the right hand side of \eqref{Eq:RicciSolitonFlow}. We denote the metrics we associate to integral curves \(\gamma:(a,b)\to\R^4\) of \(F\) (see Proposition~\ref{Prop:CurvetoMetric}) by replacing \(\gamma\) with \(g\) (e.g. if we label an integral curve \(\gamma^S\), then the associated metric is denoted \(g^S\)).
\end{notation}

 Hence, the equation we study is \begin{equation}
\gamma'= F(\gamma).
\end{equation} 

Note that if we know \(\gamma=(x,y,z,w):\R\to \R^4\), then we can reconstruct the one-parameter family \(\{g(s)\}\), and hence the metric \(g\), up to isometry.  We do this by defining \begin{equation}\label{Eq:c_h}
c^2(s)=s_0z^{-1}(s),\qquad h(s)=\int_0^s y(t)\, dt.
\end{equation}
Note that if we replaced \(h\) with \(h(s)+c\), \(c\in \R\), then we could integrate \(e^{-cD}\) to an automorphism \(\phi\) of \(\sans{S}_0\), and the map \(\Phi:\R\times \sans{S}_0\to\R\times \sans{S}_0\), \((s,x)\mapsto (s,\phi(x))\), would then be an isometry between the two metrics we obtain. 

Hence, we have the following.

\begin{prop}\label{Prop:CurvetoMetric}
An integral curve \(\gamma=(x,y,z,w)\) of \(F\) with \(z<0\) that is defined for all \(s\in \R\) determines a complete, expanding, gradient Ricci soliton metric on \(M=\R\times \sans{S}_0\) which is \(\sans{S}_0\)-invariant.  Moreover, the soliton an integral curve determines is unique up to isometry.
\begin{proof}
It remains to show that a metric \(g=ds^2+g(s)\) with \(g(s)=c^2e^{hD*}g_0\) satisfies \eqref{Eq:Cohom1ODEs(c)}. To see this, let \(X\in T_p(\sans{S}_0\cdot p)\). The formula for the mixed Ricci curvature of a Riemannian submersion \cite[9.36]{Besse87} reduce to \[\ric(X,N) = -(\dive_{s}L)(X)=-y(s)(\dive_{s} D)(X) = 0\] where the last equality follows by using \(D=\Ric -\lambda_0I\) and the Bianchi identity. Here, \(\dive_{s}\) is the divergence operator of \((\sans{S}_0,g(s))\). 
\end{proof}
\end{prop}

\begin{remark}
By the existence and uniqueness of solutions to ODEs, given an initial condition \(\gamma_0\in \R^4\), we can always find an integral curve \(\gamma\) of \(F\) with \(\gamma(0)=\gamma_0\) defined on a maximal interval of existence.  We emphasise that the non-trivial task is to find integral curves which correspond to complete metrics.  This will be the case if and only if the integral curve is defined for all \(s\in \R\); c.f. Lemma~\ref{Lemma:Completeness} in Appendix~\ref{Appendix:ProperGroupActions}.
\end{remark}

\begin{remark}
A priori, we might have expected \eqref{Eq:RicciSolitonFlow} to consist of six equations since we have three functions (\(c,h\), and \(f\)) and \eqref{Eq:Cohom1ODEs} is a second order ODE in these functions. Recall however that the gradient of the soliton potential, \(\nabla f= f' N\), is the real variable since shifting \(f\) by a constant does not affect \eqref{Eq:RS}. A similar phenomenon occurs with the variable \(h\). Changing \(h\) by a constant corresponds to pulling back the background metric \(g_0\) by an automorphism; this pullback is again a solvsoliton metric, and the corresponding cohomogeneity one manifolds are isometric.  
\end{remark}

\begin{remark}\label{Rmk:BasePoint}
In deriving the equations \eqref{Eq:RicciSolitonFlow}, we started with a cohomogeneity one manifold \((M,g)\) with \(M/\sans{S}_0=\R\) and obtain a diffeomorphism \(M\simeq \R\times \sans{S}_0\) by choosing a geodesic \(\alpha:\R\to M\) meeting all the orbits orthogonally. The point \(\alpha(0)\in M\) is mapped by this diffeomorphism to \((0,e)\in \R\times \sans{S}_0\). If instead we had chosen \(\alpha_\tau(s)=\alpha(s+\tau)\), then we obtain would have obtained a diffeomorphism of \(M\) with \(\R\times \sans{S}_0\) that maps \(\alpha_\tau(0)=\alpha(\tau)\in M\) to \((0,e)\in \R\times \sans{S}_0\).  

This corresponds to the fact that if \(\gamma\) is an integral curve of \(F\), then shifting the \(s\)-parameter gives another integral curve of \(F\). That is, if \(\gamma,\gamma_\tau\) are integral curves of \(F\) related by \(\gamma_\tau(s)=\gamma(s+\tau)\), then they correspond to a single cohomogeneity one manifold \((M,g)\) viewed from different base points. 
\end{remark}

\begin{remark}

Different choices of \(\lambda\) give equivalent systems.  We list the equivalence here for convenience. If \(\gamma_1=(x,y,z,w)\) is an integral curve of \(F_{\lambda_1}\), then, \(\gamma_2\) defined by \[\gamma_2(s) =\bigg( \sqrt{\frac{\lambda_2}{\lambda_1}}x (\tilde{s}),\sqrt{\frac{\lambda_2}{\lambda_1}}y (\tilde{s}), \frac{\lambda_2}{\lambda_1}z(\tilde{s}),\sqrt{\frac{\lambda_2}{\lambda_1}}w(\tilde{s})\bigg),\] where \(\tilde{s}=\sqrt{\frac{\lambda_2}{\lambda_1}}s\) is an integral curve of \(F_{\lambda_2}\). 

Also, if we replaced the background metric \(g_0\) with \(\rho g_0\) for some \(\rho>0\), then the variable \(y\) becomes \(\rho y\). The remaining variables remain the same.  

\end{remark}

In this paper, we are interested in the case where \(\lambda<0\). Since the cosmological constant, \(\lambda_0\), of \((\sans{S}_0,g_0)\) is also negative, we find it convenient to assume that \(\lambda=\lambda_0\). Moreover, we will assume that \(g_0\) has been scaled so that \(\tr D=1\).

\section{Existence of Complete Solutions}\label{Section:GeneralCase}

Before proving existence of complete solutions we briefly outline our approach.  Observe that existence and uniqueness of ODE solutions immediately gives us short-time solutions; the non-trivial task is to show there are solutions of \eqref{Eq:RicciSolitonFlow} which correspond to complete metrics.  A solution of \eqref{Eq:RicciSolitonFlow} gives a complete metric if and only if it is defined for all \(s\in \R\).  

There is a one-dimensional extension of \(\sans{S}_0\) that is an Einstein solvmanifold \cite[Proposition~6.1]{LafuenteLauret14}. This gives a stationary solution to \eqref{Eq:RicciSolitonFlow} which we denote \(\gamma^S\); the linearisation of the system at this stationary point has two eigenvalues with positive real part, so there is a two-dimensional unstable manifold through \(\gamma^S\).  Let us denote this unstable manifold by \(U\). Our aim is to identify a positively invariant open subset \(\Omega\) of \(\R^4\) which intersects \(U\) non-trivially where we can better understand the integral curves of \(F\).  Any initial condition in \(\Omega\cap U\) will give an integral curve \(\gamma\) of \(F\) that exists for all \(s\le 0\); by choosing the correct initial condition we can ensure \(\gamma\) exists for all \(s\ge 0\) too.

\begin{remark}\label{Rmk:Z2Symmetry}
If \(\gamma=(x,y,z,w):(\alpha,\beta)\to \R^4\) is an integral curve of \(F\), then \[\tilde{\gamma}(s) := \big(-x(-s),-y(-s),z(-s),-w(-s)\big),\quad s\in (-\beta,-\alpha),\] is also an integral curve.  Geometrically, these integrals curves correspond to the same cohomogeneity one Ricci solitons, traversed along their normal geodesics' in opposite directions. 
\end{remark}

	\subsection{Stationary Solutions}

\begin{prop}
The the stationary solutions of \eqref{Eq:RicciSolitonFlow} (with the normalisation \(\tr D=1\), \(\lambda=\lambda_0\)) are: \[\gamma^S_{\pm} : = \bigg(\pm \frac{1}{n},\pm 1, s_0,\pm 1 \bigg)\qquad  \gamma^H_{\pm} := \pm \bigg(\sqrt{\frac{-\lambda}{n}},0,0,\sqrt{-\lambda n}\bigg) \]
\end{prop}

\begin{remark}
The stationary solutions \(\gamma^S_\pm\) correspond to the same Riemannian manifold under the symmetry of the system \eqref{Eq:RicciSolitonFlow} mentioned in Remark~\ref{Rmk:Z2Symmetry}.  Therefore, we will focus our attention on \(\gamma^S_+\), and for brevity we denote this by \(\gamma^S\).  Similarly, we denote \(\gamma^H_+\) by \(\gamma^H\). The shape operator \(L\) of \(\gamma^S\) is the soliton derivation \(D\). Geometrically, \(\gamma^S\) corresponds to an Einstein solvmanifold \((\sans{S},g^S)\) obtained from \((\sans{S}_0,g_0)\) by a one-dimensional extension \cite{LafuenteLauret14}.  The stationary solution \(\gamma^H\) corresponds to the real hyperbolic space \(\R H^{n+1}\) (with sectional curvature \(\lambda/n\)).
\end{remark}

	\subsection{Long-time solutions}

\begin{prop}\label{Prop:InvariantSet}
The set 
\[\Omega=\{\gamma=(x,y,z,w)\in \R^4: 0< y<nx, s_0<z<0 \},\] is positively invariant under the flow of \(F\).
\begin{proof}
The fact \(z<0\) is preserved is obvious since \(z=0\) is invariant; clearly \(y>0\) is preserved whenever \(z<0\).

To see the final two inequalities are preserved, set \(\varphi = -x+y/n\), so that the inequality \(y<nx\) is \(\varphi<0\). We have \begin{align*}
z'& = 2z \varphi\\
\varphi'&= \frac{-\lambda}{s_0} \big(z-s_0\big) - w\varphi.
\end{align*}
The first equation shows that, in \(\Omega\), \(z'>0\), so \(z\) is increasing towards zero.  Therefore, if an integral curve \(\gamma(s)\) of \(F\) were to leave \(\Omega\) at \(s^*
\in \R\), we would have \(\varphi(s^*)=0\). But the second equation then implies \(\varphi'(s^*)<0\).  Hence, \(\Omega\) is positively invariant.
\end{proof}
\end{prop}

\begin{remark}
Since \(z(s)\) is the scalar curvature of the orbit \(\sans{S}_0\cdot p\), which is homothetic to \((\sans{S}_0,g_0)\), \(z> 0\) is not geometrically meaningful.
Moreover, the set \(\{\varphi<0,s_0<z\}\) corresponds to solutions whose orbits are expanding for all \(s\in \R\). (By this we mean that the magnitude of the scalar curvature is decreasing.) Also, since \(L = (-\varphi)I+yD\), the shape operator \(L\) will be positive definite.
\end{remark}

\begin{prop}\label{Lemma:EVals}
The unstable subspace, \(W\), at \(\gamma^S\) is two-dimensional. Hence, there is a two-dimensional unstable manifold passing through \(\gamma^S\).  Moreover, a vector \(v\in W\) is uniquely determined by its \(x\) and \(z\) coordinates and points into \(\Omega\) if and only its \(z\) coordinate is positive.

\begin{proof} The linearisation of \(F\) at \(\gamma^S\) can be calculated with a computer, so we simply state the eigenvalues without proof.  There are two distinct eigenvalues, each of multiplicity 2. These are \[ \varepsilon_\pm = \frac{-1}{2} \pm \frac{\sqrt{8+n-8s_0}}{2\sqrt{n}}.\] Note that \(8-8s_0=-8\lambda n>0\), so \(\varepsilon_+>0\) (and, of course, \(\varepsilon_-<0\)). Since \(dF(\gamma^S)\) has two positive eigenvalues, there is a two-dimensional unstable manifold \(U\) passing through \(\gamma^S\) \cite[\S2.7]{Perko01} which is tangent at \(\gamma^S\) to the \(\varepsilon_+\)-eigenspace, \(W\).  The eigenspace \(W\) is spanned by the vectors 
\[w_0 = (b,-c,p(b,c),n b ),\qquad w_1=(-a,-na,0,n),\] where \(a,b,c>0\) are fixed real numbers whose exact values are unimportant, and \(p(b,c)=2(b(n-1)+c\tr D_0^2)>0\). These are the directions from which the Einstein solution and ``non-scaling'' solution emerge; see \S\ref{Section:EinsteinCase} and \S\ref{Section:NoScal}. 

To see that a vector \(v\in W\) is  uniquely determined by it \(x\) and \(z\) coordinates, observe that \(\{e_2,e_4,w_0,w_1\}\) forms a basis for \(\R^4\) (where \(e_i\) is the \(i\)-th standard basis vector of \(\R^4\)).  Hence, if \(v\in W\) satisfies \(v\cdot e_1=v\cdot e_3=0\), then \(v=0\).  Finally, it is clear from the explicit form of \(w_0\) and \(w_1\) that if \(v\in W\) then \[\gamma^S+\tau v\in \Omega\quad \forall \tau\in (0,\delta),\] holds for some \(\delta>0\) if and only if \(v\cdot e_3>0\).  
\end{proof}
\end{prop}

\begin{lemma}\label{Lemma:ApproachLemma}
Let \(v\in W\subset \R^4\) such that \(|v|=1\).  Then there is a unique (up to shifting time) integral curve \(\gamma(s)\) of \(F\) such that \(\gamma\to \gamma^S\) as \(s\to -\infty\) and \begin{equation}\label{Eq:Emerging}
\lim_{s\to -\infty}  \frac{\gamma'(s)}{|\gamma'(s)|}=v.
\end{equation}
If \(\gamma\) satisfies \eqref{Eq:Emerging}, then we say that \(\gamma\) emerges from \(\gamma^S\) in the direction \(v\).
\begin{proof}
Since \(\gamma^S\) is a hyperbolic stationary point, an integral curve converges to \(\gamma^S\) if and only if it is in the unstable manifold. The claim now follows easily since \(F|_U\) is \(C^1\)-equivalent in a neighbourhood of \(\gamma^S\) to the linear system \(\gamma' = \varepsilon_+ \gamma\); see the final theorem in \S2.8 of \cite{Perko01}.
\end{proof}
\end{lemma}

\begin{corollary}\label{Corollary:Uniqueness}
Let \(\gamma_1\) and \(\gamma_2\) be integral curves of \(F\) emerging from \(\gamma^S\). Then, \(z'_i>0\) for \(s<0\) and \[\lim_{s\to -\infty} \frac{x'_1(s)}{z'_1(s)}=\lim_{s\to -\infty} \frac{x'_2(s)}{z'_2(s)}, \] implies that \(\gamma_1(s)=\gamma_2(s+\tau)\) for some \(\tau\in \R\).
\begin{proof}
By Proposition~\ref{Lemma:EVals}, \(W\) can be parametrised by the \(xz\)-coordinates. The assumptions on \(z'_i\) and \(x'_i\) uniquely determine a unit vector in the \(xz\)-plane, and therefore also uniquely determine a unit vector \(v\in W\). Hence, \(\gamma_1\) and \(\gamma_2\) must emerge from \(\gamma^S\) in the same direction, so \(\gamma_1(s)=\gamma_2(s+\tau)\) for some \(\tau\in \R\) by Lemma~\ref{Lemma:ApproachLemma}.
\end{proof}
\end{corollary}

If we parametrise the unit circle in the \(xz\)-plane by the angle \(\theta\in [-\pi,\pi)\) that \(v\) makes with the positive \(z\)-axis, then the vectors that point into \(\Omega\) correspond to \(\theta \in (\frac{-\pi}{2},\frac{\pi}{2})\).  Moreover, the vectors \(w_0\) and \(w_1\) correspond to \(\theta=\theta_0>0\) and \(\theta= -\frac{\pi}{2}\) respectively. (\(\theta_0\) could be calculated explicitly, but its exact value is unimportant.) It remains to show which directions give complete solutions.

\begin{notation}
We denote by \(\gamma_\theta:(-\infty,s^*)\to \R^4\), \(s^*\in (-\infty,\infty]\), an integral curve of \(F\) emerging from \(\gamma^S\) in the direction \(\theta\), defined on its maximal interval of existence.
\end{notation}

\begin{theorem}\label{Thm:MainExistence}
Let \(\theta\in (-\frac{\pi}{2},\theta_0)\), let \(v_\theta\in W\) the vector corresponding to \(\theta\), and let \(\gamma_\theta=(x,y,z,w)\) be an integral curve of \(F\) emerging from \(\gamma^S\) in the direction \(v_\theta\).  Then, \(\gamma_\theta(s)\in \Omega\) and \(w(x)>x(s)n\) for all \(s\in (-\infty,s^*)\), where \(s^*\) is maximal existence time. Moreover, it follows that \(s^*=\infty\).

\(\gamma_\theta\) is defined for all \(s\in \R\).  Moreover, \(\gamma_\theta(s)\in \Omega\) and \(w(x)>x(s)n\) for all \(s\in \R\).

\begin{proof}
Since \(\Omega\) is defined by linear inequalities, and the trace of\(\gamma_\theta\) is well approximated by the line through \(\gamma^S\) in the direction \(v_\theta\) as \(s\to -\infty\), it follows easily that \(\gamma_\theta(s)\in \Omega\) for sufficiently negative \(s\). That \(\gamma_\theta(s)\in \Omega\) for all \(s\in (-\infty,s^*)\) now follows immediately since \(\Omega\) is positively invariant. 

Note that \(v_\theta\) points into the region \(\{xn<w\}\cap \Omega\) precisely when \(\theta\in (-\frac{\pi}{2},\theta_0)\), so the argument above also shows that there exists \(s_1\in \R\) such that \(x(s)n<w(s)\) for all \(s<s_1\). It remains to show this continues to hold.

Since \(\gamma_\theta\to \gamma^S\) exponentially as \(s\to -\infty\), \(\varphi(s):=w(s)-x(s)n\) converges to \(0\) exponentially, and hence is integrable on \((-\infty,s]\).  Let \[\Phi(s)=\int_{-\infty}^s \varphi(\tau)\, d\tau.\]  Observe that, unwinding the definitions, \(\Phi\) is the soliton potential normalised in such a way that \(\Phi\to 0\) as \(s\to -\infty\).  In particular, \(\Phi\) satisfies the following equation (see, for example, Equation 3.13 in \cite{DancerWang10}) \begin{equation}\label{Eq:Cohom1Bianchi}
\Phi''+w\Phi'+2\lambda\Phi=C, \quad C\in \R.
\end{equation} (Note that in \cite{DancerWang10}, \eqref{Eq:Cohom1Bianchi} is derived from \eqref{Eq:Cohom1ODEs}, so it does not require completeness.)  Since the left-hand side converges to \(0\) as \(s\to -\infty\), we must have that \(C=0\). It is easy to see from \eqref{Eq:Cohom1Bianchi} that the set \(\{\Phi>0,\Phi'>0\}\) is preserved; it is also easy to see that \(\Phi(s_1)\), \(\Phi'(s_1)>0\) for \(s_1\in \R\) chosen so that \(x(s)n<w(s)\) for all \(s<s_1\). Hence, the claim holds since \(0<\Phi'(s)=x(s)n-w(s)\) for all \(s\) such that \(\gamma_\theta\) is defined.

Finally, we must have \(s^*=\infty\) since \(w\) grows at most linearly, \(0<y<nx<w\), and \(s_0<z<0\).
\end{proof}
\end{theorem}

Note that once we have \(\Phi'(s)>0\) for all \(s\in \R\) we can use \eqref{Eq:Cohom1Bianchi} and argue as in \cite[Proposition~1.11]{BuzanoDancerWang15II} to show \(\Phi''(s)>0\) for all \(s\in \R\). Since \(L\) is also positive definite, this implies that the solitons given by Theorem~\ref{Thm:MainExistence} will have a convex potential function. This agrees with \cite{MaschlerReam22}. (Note that \cite{MaschlerReam22} uses a different sign convention, so the soliton potential is concave for them.)

\begin{corollary}
The Hessian of the soliton potential is positive definite for the solitons produced in Theorem~\ref{Thm:MainExistence}. Hence, the soliton potential is convex.
\end{corollary}

	\subsection{Asymptotics}

Now that we have existence of complete solitons, it is interesting to ask how the metrics look as the parameter \(s\) approaches \(\pm \infty\).

Recall that \(\gamma^S\) is the stationary solution \(\gamma^S=(n^{-1},1,s_0,1)\). Hence, as \(s\to -\infty\), we have \(c=s_0 z^{-1}\to 1\). Moreover, if \(h\) is an anti-derivative of \(y\), then \(y\to 1\) implies \(h\sim s\) as \(s\to -\infty\). Therefore, if \(g\) is a metric associated to a curve \(\gamma\) that is in the unstable manifold at \(\gamma^S\) then \(g\) is asymptotic at the corresponding end to the Einstein solvmanifold metric \[g^S=ds^2+e^{sD_*}g_0.\]

\begin{remark}\label{Rmk:FormalAsymptotics}
If we consider the curves \(\gamma_\tau(s)=\gamma(s+\tau)\), then we think of the metrics \(g_\tau=ds^2+c_\tau e^{h_\tau D_*}g_0\) defined by \(\gamma_\tau\) as being the same metric viewed from different base points (see Remark~\ref{Rmk:BasePoint}). Since \(c(s)=s_0 z(s)^{-1}\to 1\) as \(s\to -\infty\), \(c_\tau(s)\to 1\) for all \(s\in \R\) as \(\tau\to-\infty\). Similarly, \(y\to 1\) implies \(h'_\tau(s)=y_\tau(s)\to 1\) for all \(s\), and hence \(h_\tau(s)\to s\).

Therefore, we can view the asymptotics above in the following way: if we take a sequence of points \(\{p_i\}\) in \((M,g)\) that converge to the Einstein solvmanifold end of \((M,g)\), then the sequence of pointed manifolds \((M,g,p_i)\) converges in the pointed Cheeger-Gromov topology to \((\sans{S},g^S,e)\).
\end{remark}

We now turn our attention to computing the asymptotics at the other end of \(M\).  We start by finding the limits of the quantities \(x,y,z,w\).

\begin{prop}
Let \(\gamma_\theta=(x,y,z,w)\) be the integral curve of \(F\) emerging from \(\gamma^S\) in the direction \(\theta\in (-\frac{\pi}{2},\theta_0)\).  Then, there exists \(C>0\) such that \(x,y\le C\) for all \(s\in \R\).  Moreover, we have the following limits as \(s\to \infty\):
\begin{equation*}
 x,y\to 0,\qquad
w\to +\infty .
\end{equation*}
\begin{proof}
We first show that \(x,y\) are bounded; observe that it suffices to show that \(x\) is bounded. But if \(x^2n>-\lambda\) then using that \(nx<w\) we have \[x'=\frac{z}{n}-\lambda -wx< -\lambda -nx^2<0.\] Hence, \(x\) is bounded.

Next, we show that \(w\) diverges to \(\infty\). Since \((M,g)\) is complete, \(|\nabla f|=f'\) must be unbounded \cite[Theorem~2]{PigaloRimoldiSetti11}. (Note that \(\theta<\theta_0\) ensures the curve \(\gamma\) does not correspond to an Einstein metric.) This forces \(w=nx+f'\) to be unbounded since \(x>0\). By differentiating \eqref{Eq:RicciSolitonFlow(d)} and using \eqref{Eq:RicciSolitonFlow(a)}, \eqref{Eq:RicciSolitonFlow(b)} we have \[w'' = -2xx'n-2yy' =-2x(z-\lambda n-wxn)-2y(zs_0^{-1}-wy)\tr D_0^2 >2w(x^2n +y^2\tr D_0^2) -C,\] for some \(C>0\) since \(x,y,z\) are bounded. Hence, if \(w>(-2\lambda)^{-1}C\) then \(w'>0\) is preserved, so \(w\) is eventually monotone increasing. Hence, \(w\) is unbounded, and is eventually monotone increasing, so it must diverge to \(+\infty\).

Finally, we show that \(x,y\to 0\). Assume \(x\ge \varepsilon>0\). By choosing \(s^*\) large enough so that \(w\) is monotone increasing and \(w>-\lambda/\varepsilon+1\), we find that \[x' =\frac{z}{n}-\lambda -wx < -\varepsilon,\qquad \forall s>s^*. \] But this clearly gives a contradiction.  Hence, \(x\) converges to \(0\), which forces \(y\to 0\) as \(s\to \infty\).
\end{proof}
\end{prop}

Observe that since \(z<0\) and \(z'>0\), it is clear that \(z\to z_0\le 0\) as \(s\to \infty\).  \\

Since the variable \(w\) is diverging to \(+\infty\), it will be useful to introduce the variable \(v:=w^{-1}\) and to take a new time parameter \(\tau\) defined by \(\tau'=w\).  We can then study the behaviour of \((x,y,z,v)\) as \(\tau \to \infty\) using centre manifold theory.  We state the results here (translated back into our original variables), but defer the calculation to Appendix~\ref{Appendix:CentreMfldCalcs}.

\begin{prop}\label{Prop:InftyAsymptotics}
If \(\gamma_\theta=(x,y,z,w):\R\to \R^4\), \(\theta\in (0,1)\), then, \(z\to 0\) as \(s\to \infty\) and
\begin{equation*}
w\sim -\lambda s,\quad
x \sim s^{-1},\quad
y\cdot s^{2}\to 0,\quad
z\sim -\alpha s^{-2},\quad \alpha>0,
\end{equation*}
as \(s\to \infty\). 
\end{prop}
Therefore, the metric \(g_\theta\) associated \(\gamma_\theta\) is asymptotic to the cone metric \[g_c = ds^2+\alpha|s_0|s^2 g_0, \qquad \alpha>0.\]

\begin{remark}
Similar to Remark~\ref{Rmk:FormalAsymptotics} above, it is not difficult to formalise the previous asymptotics. More precisely, if we identify \(M\) with \(\R\times \sans{S}_0\) by picking a geodesic meeting the orbits orthogonally (and pick the direction we are travelling along this geodesic to match our convention above), then
\[\frac{1}{\tau^2}\rho_\tau^*\Phi^*g \to ds^2+\alpha|s_0|s^2g_0,\qquad\tau \to \infty,\] 
where \(\rho_\tau:(0,\infty)\times \sans{S}_0 \to (0,\infty)\times \sans{S}_0\) is defined by \(\rho_\tau(s,x) = (\tau s, x)\) for \(\tau>0\). Here \(\Phi(s,x)=(s,\phi(x))\), \(d\phi = e^{h_\infty D}\), where \(h_\infty\in \R\) is the limit of \(h\) as \(s\to \infty\) (which exists \(y\cdot s^2\to 0\)). 

If \(\sans{S}_0\) admits a cocompact lattice, \(\Gamma\), then \( M/\Gamma \approx_{\text{diffeo.}} \R\times  (\sans{S_0}/\Gamma)\) has two non-compact ends. In this case, \(M/\Gamma\) is asymptotic to a cone metric over \(\sans{S}_0/\Gamma\) along one of these ends in the sense of \cite[Definition~1.1]{KotschwarWang15}.
\end{remark}

\begin{theorem}
Let \(\theta_1,\theta_2\in (-\frac{\pi}{2},\theta_0)\), and let \(\gamma_1,\gamma_2:\R\to \R^4\) be integral curves of \(F\) emerging from \(\gamma^S\) in the direction \(v_{\theta_1},v_{\theta_2}\). If the Riemannian manifolds \((M_1,g_1)\) and \((M_2,g_2)\) determined by \(\gamma_1\) and \(\gamma_2\) are isometric, then \(\gamma_1(s)=\gamma_2(s+\tau)\) for some \(\tau\in \R\).
\begin{proof}
Let \(\phi:M_1\to M_2\) be an isometry. Observe that in the set-up in \S2 any unit-speed geodesic in \((M,g)\) meeting all the orbits would give rise to the same integral curve of \(F\), up to shifting or reversing the \(s\)-parameter. And, if \(\alpha_1\) is a geodesic in \(M_1\) meeting the \(\sans{S}_0\)-orbits orthogonally, then \(\phi\circ\alpha_1\) is a geodesic in \(M_2\) meeting the \(\sans{S}_0\)-orbits orthogonally; see Proposition~\ref{Prop:OrbitPreservation}. Hence, we can assume the geodesics, \(\alpha_i\), we picked for \((M_i,g_i)\), \(i=1,2\), are related by \(\alpha_2(s+\tau)=\phi(\alpha_1(s))\) for some \(\tau\in \R\). (The \(s\)-parameter will not be reversed since the asymtotics must match).

Since \(z_i\) and \(x_i\) are the scalar curvature and mean curvature of the orbits restricted to the geodesic \(\alpha_i\), it follows that \(z_1(s)=z_2(s+\tau)\) and \(x_1(s)=x_2(s+\tau)\). Hence, 
\[\lim_{s\to -\infty} \frac{x_1(s)}{z_1(s)}=\lim_{s\to -\infty}\frac{x_2(s)}{z_2(s)},\] so Corollary~\ref{Corollary:Uniqueness} implies \(\gamma_1(s)=\gamma_2(s+\tau)\) for some \(\tau\in \R\).
\end{proof}
\end{theorem}

\section{The Einstein Case}\label{Section:EinsteinCase}

We turn our attention now to the existence of solutions that define Einstein manifolds. Our main aim in this section is to show that there is a cohomogeneity one Einstein manifold interpolating between the Einstein solvmanifold \((\sans{S},g^S)\) and real hyperbolic space \((\R H^{n+1},g^H)\). In the case that \(\sans{S}_0=H_3(\R)\) is the Heisenberg group, this is the one-loop deformed universal hypermultiplet metric studied in \cite{CortesSaha22}. In this case, \(\sans{S}=\C H^2\).

Looking for Einstein metrics amounts to finding integral curves of \(F\) in the invariant subset \[E=\{(x,y,c,w): xn=w, z = \lambda(n-1)+x^2n(n-1)-y^2\tr D_0^2\}.\] 
The condition \(w=xn\) simply states that \(f'\equiv 0\), while the second condition is a conservation law that is always true for cohomogeneity one Einstein manifolds \cite[Remark~2.2]{EschenburgWang00}.

 The vector field \(F|_{E}\) becomes \begin{align*}
x'&=-\frac{\lambda}{n}-x^2-y^2\frac{\tr D_0^2}{n}=:\mathcal{E}_1(x,y)\\
y'&=\frac{\lambda(n-1)}{s_0} +x^2\frac{n(n-1)}{s_0}+\frac{y^2}{n}- nxy=:\mathcal{E}_2(x,y).
\end{align*}

	\subsection{Existence of an Einstein Solution}
	
\begin{theorem}\label{Thm:EinsteinSoln}
There is an integral curve of \(\mathcal{E}\) that joins the stationary solution \(\gamma^{S}\) to the stationary solution \(\gamma^H\).
\begin{proof}
Consider the set \[K=\big\{(x,y)\in \R^2:x,y\ge 0,\mathcal{E}_1(x,y)\ge 0 ,\mathcal{E}_2(x,y)\le 0 \big\}.\]

The matrix \(d\mathcal{E}(\gamma^S)\) has eigenvalues \(\varepsilon_\pm\), both with multiplicity one.  The vector \[w_0=-\bigg( n-4-n\sqrt{-8\lambda n+1},\frac{4n(n-1)}{-s_0}+2n^2\bigg),\] is an eigenvector corresponding to \(\varepsilon_+\). Since \(\sqrt{-8\lambda n+1}>1\), \(w_0\) points into the region \(K\) and hence the unstable manifold at \(\gamma^S\) intersect \(K\). (In particular, \(K\) is non-empty.)

The set \(K\) is also positively invariant: if \((x,y)\) is a point in the boundary of \(K\) other than \(\gamma^S\) or \(\gamma^H\), then either \(\mathcal{E}_1(x,y)=0\) or \(\mathcal{E}_2(x,y)=0\). The derivative of \(\mathcal{E}_1\) (resp. \(\mathcal{E}_2\)) in the direction of \(\mathcal{E}(x,y)\) is then positive (resp. negative), so integral curves cannot leave \(K\).

Clearly \(K\) is compact, so if \(\gamma\) is an integral curve whose image is the part of unstable manifold that is in \(K\), then \(\gamma\) is defined for all \(s\in \R\).

We have \(\gamma\to \gamma^S\) as \(s\to -\infty\) since it is in the unstable manifold. To see that \(\gamma\to\gamma^H\) as \(s\to \infty\), observe that \(K\) is precisely the set in which \(x'\ge 0\) and \(y'\le 0\). By monotonicity, \(x\) and \(y\) must converge; continuity implies \(x'(s)=\mathcal{E}_1(x(s),y(s))\) and \(y'(s)=\mathcal{E}_2(x(s),y(s))\) also converge. The limit of \(x'\) and \(y'\) must be \(0\), so we conclude that \(\gamma\) converges to a stationary solution. The only possibilities are \(\gamma^S\) and \(\gamma^H\), and it clearly cannot be \(\gamma^S\).
\end{proof}
\end{theorem}

\begin{remark}
Recall that the Damek-Ricci spaces are non-symmetric, harmonic manifolds. They are also Einstein solvmanifolds, so Theorem~\ref{Thm:EinsteinSoln} gives inhomogeneous deformations of the Damek-Ricci spaces. Unfortunately, these deformations are not harmonic: the quantity \(|R_N|^2\), where \(R_vw=R(v,w)v\) and \(N_s\) is the unit normal field of \(\{s\}\times \sans{S}_0\), is not constant. 
\end{remark}

	\subsection{Asymptotics of Einstein solutions}

Since the integral curve \(\gamma=(x,y):\R\to \R^2\) from Theorem~\ref{Thm:EinsteinSoln} joins \(\gamma^S\) to \(\gamma^H\), we expect the geometry of the manifold \((M,g)\) determined by \(\gamma\) to approach the geometry of \((\sans{S},g^S)\) and \(\R H^{n+1}\) as \(s\to -\infty\) and \(s\to \infty\) respectively. As \(s\to -\infty\), this follows from the previous section. To see \(\R H^{n+1}\) as \(s\to \infty\), we need to account for the fact the orbits are expanding. We deal with this by looking in the correct coordinates.

As before, for \(\tau \in \R\) let \(\gamma_\tau=(x_\tau,y_\tau):\R\to \R^2\) be the curve \(\gamma_{\tau}(s):=\gamma(s+\tau)\). Recall that each \(\gamma_\tau\) gives rise to the same abstract manifold, \((M,g)\), viewed from  different base points.

Consider the following coordinates on \(\R\times \sans{S}_0\). The Lie algebra is a semi-direct product \(\mathfrak{s}_0=\mathfrak{a}\oplus \mathfrak{n}\), where \(\mathfrak{n}\) is the nilradical and \(\mathfrak{a}\) is the orthogonal complement of \(\sans{s_0}\) with respect to \(g_0\). The subalgebra \(\mathfrak{a}\) is in fact abelian. The map \(\exp:\mathfrak{s}_0=\mathfrak{a}\oplus \mathfrak{n}\to \sans{S}_0\), \[(A,X)\mapsto \exp A \exp X,\] is a diffeomorphism onto \(\sans{S}_0\).

Let \(\{A_i^0\}\cup\{X_i^0\}\) be a \(g_0\)-orthonormal basis of \(\mathfrak{s}_0\), and let \(X_i^\tau:=c^{-1}(\tau)X_i^0\), \(A_i^\tau=c^{-1}(\tau)A_i^0\). The map \(\Phi_\tau:\R^{n+1}\to \R\times \sans{S}_0\), where \[\Phi_\tau(s,x) = \big(s,\exp(x^i A_i^\tau)\exp( x^jX_j^\tau)\big), \] is a diffeomorphism that takes \(0\in \R^{n+1}\) to \((0,e)\in \R\times \sans{S}_0\). The pullback of \(g^\tau\) by \(\Phi_\tau\) is \[\Phi_\tau^*g^\tau = ds^2 +\frac{c^2_\tau(s)}{c^2_\tau(0)}e^{h_\tau(s)D*}\bigg(g_{ij}(x)dx^idx^j\bigg),\]
where \(g_{ij}(x)\) is an analytic function of \(x\) whose coefficients are given by a universal expression in the structure constants \(\{\mu_{ij}^k\}\) of the basis \(\{A_i^\tau\}\cup \{X_i^\tau\}\); see \cite[Proposition~6.9]{Lauret12}. The \(c^2_\tau(0)\) term in the denominator appears since \cite{Lauret12} uses a metric making the basis of \(\mathfrak{s}_0\) orthonormal, so we have applied it to \(c^2_\tau(0)g_0\). In particular, it follows that coefficients \(g_{ij}(x)\), and their derivatives, converge to the coefficients of the Euclidean metric, and their derivatives, uniformly on compact subsets as the structure constants converge to \(0\). The structure constants of \(\{A_i^\tau\}\cup \{X_i^\tau\}\) are related to the structure constants of \(\{A_i^0\}\cup\{X_i^0\}\) by a factor of \(c^{-1}(\tau)\), so they go to \(0\) as \(\tau\to \infty\).

Now, it is not difficult to check that \(h_\tau\) converges to \(0\) uniformly on any compact subset of \(\R\). Moreover, since \(c' = (x-y/n)c\), we deduce that \[c(s)\sim \rho_0 e^{s\alpha},\quad s\to \infty,\] for some \(\rho_0>0\) and \(\alpha:=\sqrt{\frac{-\lambda}{n}}\). Thus, \[\frac{c_\tau(s)}{c_\tau(0)}=\frac{c(s+\tau)}{c(\tau)}\to e^{s\alpha},\quad \tau\to \infty.\]
In summary, we have the following.

\begin{prop}
With the notation above, 
\[\Phi_\tau^*g^\tau \to g^H = ds^2 +e^{2s\alpha}\big( (dx^1)^2+\ldots +(dx^n)^2\big),\] as \(\tau \to \infty\).
This is the real hyperbolic metric of sectional curvature \(-\alpha^2 = \frac{\lambda}{n}<0\).
\end{prop}

\section{Ricci solitons whose orbits are isometric to Solvsolitons}\label{Section:NoScal}

There is a second invariant of subset of the vector field \(F\). This subset corresponds to the ansatz that the shape operator \(L\) is a multiple of the soliton derivation \(D\) for all \(s\in \R\).  In terms of the geometry of the orbits, this means that we require our orbits to be isometric to the background solvsoliton \((\sans{S}_0,g_0)\) (i.e. \(g(s)=e^{h(s)D*}g_0\)). In the variables \((x,y,z,w)\), this subset is given by \[\{(x,y,z,w): y\tr D=xn, z=s_0\},\] where \(s_0=\lambda n+\tr D\) is the scalar curvature of \((\sans{S}_0,g_0)\); it is easy to see from \eqref{Eq:RicciSolitonFlow} that this subset is invariant.

The vector field \(F\) restricted to this subset gives the system \begin{subequations}\label{Eq:NoScalSystem} \begin{align}
\label{Eq:NoScalSystem(a)} y'& = 1- wy\\
\label{Eq:NoScalSystem(b) }w'&= -\lambda(1 - y^2)
\end{align}\end{subequations}

Heretofore, we have assumed that \((\sans{S}_0,g_0)\) is not scalar flat; this is because the flat case does not fit into our framework. A scalar flat solvsoliton is Ricci flat. In this case, the only way to have a non-trivial soliton derivation is if \(D=-\lambda I\), \(\lambda<0\). But \(D=-\lambda I\) is a derivation if and only if \(\sans{S}_0=\R^n\) is abelian. This reduces our ansatz to a warped product ansatz. For \(\sans{S}_0=\R^2\), this was studied in \cite{Ramos13}. Ramos shows there is a unique integral curve of \eqref{Eq:NoScalSystem} corresponding to a complete metric with bounded sectional curvature. In fact, Ramos shows any non-flat Ricci soliton on \(\R\times T^2\) whose curvature is bounded from below must a quotient of the soliton in \cite{Ramos13}.

Since \cite{Ramos13} studies the system \eqref{Eq:NoScalSystem} using phase plane analysis for \(\lambda=-1/2\), and it is easy to see that Ramos' analysis remains valid for \(0<-\lambda<1\), we state our results here without proof. Note that Ramos uses variables, \(F\) and \(H\), that are related to ours by \(y=2H\), \(w=2H-F\). (Note that \(-1<\lambda<0\) does hold, since the eigenvalue of \(D\) are non-negative and \(\tr D=1\).)

\begin{prop}\label{Prop:NOScalExistence}
The system\begin{align*}
y'& = 1- wy\\
w'&= -\lambda(1 - y^2)
\end{align*}  
admits an integral curve, \(\gamma_1\), defined for all \(s\in \R\).  This curve is contained in global unstable manifold emanating from the stationary solution \((1,1)\in \R^2\), and has the following asymptotics: \[w\sim -\lambda s,\qquad y\sim \frac{1}{-\lambda s},\qquad s\to \infty.\]
Hence, if \((M,g)\) is the Riemannian manifold induced by \(\gamma_1\), then \(g\) is asymptotic as \(s\to \infty\) to the metric
\[ ds^2+e^{\frac{\log s}{-\lambda}D_*}g_0.\] 
\end{prop}

The integral curve \(\gamma\) given by Proposition~\ref{Prop:NOScalExistence} emerges from \(\gamma^S\) in the direction \(w_1=(-a,-na,0,n)\) where \(a>0\) is a constant. Observe that when \(\sans{S}_0=\R^n\), \(D=-\lambda I\), so these metrics are asymptotically conical in this case.

\appendix
\section{Centre Manifold Calculations}\label{Appendix:CentreMfldCalcs}

We give in this appendix the calculations which were omitted from Section~\ref{Section:GeneralCase}. We use the theory of centre manifolds; see \cite{Carr82,Bressan07}. Recall from \S\ref{Section:GeneralCase} that the integral curves of \(F\) that we are interested in approach \((0,0,z_0,\infty)\) as \(s\to \infty\). By defining new variables (see the next paragraph) this becomes a non-hyperbolic stationary point with two negative eigenvalues and two zero eigenvalues. Therefore, there is a two-dimensional invariant centre manifold. Any integral curve converging to the stationary point that is not in the stable manifold converges exponentially to a trajectory in the center manifold; this means that understanding the trajectories in the center manifold is enough to determine the asymptotics.

Consider the system 
\begin{subequations}\label{Eq:NewRicciSolitonFlow}
\begin{align}
\label{Eq:NewRicciSolitonFlow(a)} \frac{d}{d\tau} x &= \frac{zv}{n}-\lambda v-x,\\
\label{Eq:NewRicciSolitonFlow(b)} \frac{d}{d\tau}y&=\frac{zv}{s_0}-y,\\
\label{Eq:NewRicciSolitonFlow(c)} \frac{d}{d\tau}z& = 2vz\bigg( \frac{y}{n}-x\bigg),\\
\label{Eq:NewRicciSolitonFlow(d)} \frac{d}{d\tau}v& =-v^3\bigg( -\lambda -nx^2- y^2\tr D_0^2\bigg).
\end{align}
\end{subequations}

This is simply system \eqref{Eq:RicciSolitonFlow} with the variable \(w\) replaced by \(v=w^{-1}\) and reparametrised with the variable \(\tau\) defined by \(\tau'=v^{-1}=w\).  We are interested in the dynamics of integral curves as they approach stationary solutions of the form \(\gamma_\infty=(0,0,z_0,0)\), \(z_0\le 0\).  

The linearisation of \eqref{Eq:NewRicciSolitonFlow} at \(\gamma_\infty\) has an eigenvalue of \(-1\) with multiplicity 2 and an eigenvalue of \(0\) with multiplicity 2.  The eigenspaces are \[V_0 = \R\{ae_1+be_2+e_4,e_3\},\quad V_{-1}=\R\{e_1,e_2\},\] where \(e_i\in \R^4\) is the \(i\)-th standard basis vector and \[a:=\frac{z_0}{n}-\lambda >0, \qquad b:= \frac{z_0}{s_0}\ge 0.\]  Note that if \(z_0=0\), then \(a=-\lambda\) and \(b=0\).

The stable manifold at \(\gamma_\infty=(0,0,z_0,0)\) is given by \[S =\{((x,y,z_0,0):x,y\in \R\}. \] (This subset in invariant and explicitly solving the restricted system shows it must be the stable manifold by uniqueness.)  Therefore, the behaviour of integral curves \((x,y,z,v)\) of \eqref{Eq:NewRicciSolitonFlow} that converge to \((0,0,z_0,0)\) where \(v>0\), \(z<0\) will be determined by the centre manifold.

Changing variables to \((\xi,\eta)\in \R^2\times \R^2\) where \[\eta_1 = x-av,\ \eta_2= y- bv,\ \xi_1=z-z_0,\ \xi_2=v,\] puts the system into normal form:
 \begin{subequations}\label{Eq:NewRicciSolitonFlow_NormalForm}
 \begin{align}
\frac{d}{d\tau}\xi&= F(\xi ,\eta )\\
\frac{d}{d\tau}\eta &= - \eta + G(\xi ,\eta).
\end{align}\\
\end{subequations}
Here \(F,G:\R^2\times \R^2\to \R^2\) are given by \begin{align*}
F(\xi,\eta)& =\big(2(\xi_1+z_0)\xi_2({\textstyle\frac{\tr D}{n}}(\eta_2+b\xi_2)-\eta_1-a\xi_2), P(\xi,\eta)\big) \\
G(\xi,\eta)& =   \big(\xi_1\xi_2 + aP(\xi,\eta),{\textstyle\frac{1}{s_0}}\xi_1\xi_2-{\textstyle\frac{z_0}{s_0}}P(\xi,\eta)\big), 
\end{align*}
and the function \(P:\R^2\times \R^2\to \R\) is defined by
\[P(\xi,\eta) = -\xi_2^3\bigg( -\lambda - n\big(\eta_1+a\xi_2\big)^2 - \big(\eta_2+b\xi_2\big)^2\tr D_0^2\bigg).\]
Centre manifold theory ensures that there is a \(C^2\) function \({h}:\R^2\to \R^2\) whose graph is the centre manifold on a neigbourhood of \(0\in \R^2\) \cite[Theorem~2]{Carr82}. This \({h}\) is a solution to the PDE \begin{equation}\label{Eq:centremanifoldPDE}
d{h}_\xi[ F(\xi,{h}(\xi))] +{h}(\xi) -G(\xi,{h}(\xi))=0,
\end{equation} with \({h}(0)=0,dh_0=0\), where \(d{h}_\xi:\R^2\to \R^2\) is the Jacobian matrix of \({h}\) at \(\xi\).

Of course, solving \eqref{Eq:centremanifoldPDE} for \(h\) explicitly is not possible. However, for our purposes the following information about \(h\) will be sufficient.

\begin{lemma}
There is a \(C^1\) function \(g:\R^2\to \R^2\) such that \(g(0)=0\) and \[{h}(\xi)= \xi_2 g(\xi),\] in a neighbourhood of the origin.
\begin{proof}
This follows by rearranging \eqref{Eq:centremanifoldPDE} for \({h}\) (since we know \({h}\) exists and satisfies \eqref{Eq:centremanifoldPDE}): \[ {h}(\xi) =G(\xi,{h}(\xi))-d{h}_\xi[ F(\xi,{h}(\xi))] = \xi_2 g(\xi).\]  
The function \(g\) is continuously differentiable at \(\xi=0\) since \(G\) and \(F\) both have a factor of \(\xi_2\) (i.e. \(G = \xi_2 \tilde{G}\)) and all the terms on the right hand side are \(C^1\). And \(g(0)=0\) follows by substitution.
\end{proof}
\end{lemma}

The flow on the centre manifold is defined by the system \cite[Section~2.12]{Perko01} \[\frac{d}{d\tau}\xi = F(\xi,h(\xi)). \] 
Observe that \(\eta/\xi_2\to 0\) as \(\xi\to 0\) in the centre manifold since \(\eta =h(\xi) = \xi_2g(\xi)\). Therefore, if \((\xi(\tau),\eta(\tau))\) is an integral curve converging to \(0\) as \(\tau\to \infty\) that is not in the stable manifold, then by dividing \(\eta_1,\eta_2\) by \(\xi_2\), we obtain that \(x\sim av\) and \(y\sim bv\) (when \(z_0\neq 0\)).

Up until now, we have only assumed that \(z_0\le 0\).  We show now that for the integral curves we are interested in, the only possibility is \(z_0=0\).

\begin{lemma}
Let \(\gamma=(x,y,z,v)\) be an integral curve of \eqref{Eq:NewRicciSolitonFlow} converging to \((0,0,z_0,0)\) as \(\tau\to \infty\).  Assume further that \(\gamma\) is not is the stable manifold at \((0,0,z_0,0)\).  Then, either \(z_0=0\), or \(z_0=s_0\). 
\begin{proof}
Suppose for a contradiction that there is \(z_0\neq 0,s_0\), and \(\gamma=(x,y,z,w)\) such that \(\gamma= (x,y,z,w)\to (0,0,z_0,0),\) but \(\gamma\) is not in the stable manifold at \((0,0,z_0,0)\). Let \(\phi=(\xi,h(\xi))\) be a trajectory in the centre manifold such that \(|\gamma(\tau)-\phi(\tau)|=O(e^{-\rho \tau})\) as \(\tau\to \infty\) where \(\rho>0\). We obtain, in particular, that \[x\sim av,\qquad y\sim bv,\qquad v\sim (-2\lambda\tau)^{\frac{-1}{2}},\quad \tau\to \infty.\]
Hence, \[\frac{\frac{d}{d\tau}z}{z} \sim 2v^2\bigg( \frac{\tr D}{n}b-a\bigg) \sim \frac{1}{-\lambda \tau} \bigg(\frac{z_0}{s_0}-1\bigg)= \frac{-C}{\tau}, \ C>0.\]
We now obtain a contradiction from the limit comparison test for integrals: \[\int_1^\infty \frac{\frac{d}{d\tau}z}{z}  \,d\tau,\quad \text{and}\quad \int_1^\infty \frac{1}{\tau}\, d\tau,\] should converge or diverge together.
\end{proof}
\end{lemma}
Recall that if \(z_0=0\), then the constants \(a,b\) appearing in the linearisation of \eqref{Eq:NewRicciSolitonFlow} are \(a=-\lambda\) and \(b=0\) respectively. In this case, the flow on the centre manifold takes the form \begin{subequations}\label{Eq:CentreManifoldSystem} \begin{align}
\xi_1' &= 2\xi_1\xi_2^2(\lambda + r_1(\xi)),\\
\xi_2'&= -\xi_2^3(-\lambda +r_2(\xi)),
\end{align}
\end{subequations}
where \(r_1,r_2\) are both \(C^1\) and \(r_1(0)=r_2(0)=0\).
\begin{lemma}\label{Lemma:AsymptoticsofScalarCurv}
Let \(\xi\) be a solution of this system, converging to \(0\in \R^2\). Assume further that the curve \(\xi\) is contained in the quadrant \(\{\xi_1<0,\xi_2>0\}\). Then, \[\lim_{\tau \to \infty} \frac{\xi_1}{\xi_2^2}=-\alpha,\] for some \(\alpha>0\).
\begin{proof}
If we reparametrise by \(t\), where \(\xi_2^2 d\tau = dt\), the system \eqref{Eq:CentreManifoldSystem} has a hyperbolic stationary point at \(0\). The linearisation at this point is block diagonal with eigenvalues \(2\lambda\) (in the \(\xi_1\) direction) and \(\lambda\) (in the \(\xi_2\) direction). Hence, the claim follows from the Hartman-Grobman theorem \cite[Section~2.12]{Perko01}.
\end{proof}
\end{lemma}

\begin{prop}
Let \(\gamma=(x,y,z,v)\) be an integral curve of \eqref{Eq:NewRicciSolitonFlow} converging to \((0,0,0,0)\) as \(\tau\to \infty\) that is not in the stable manifold. Then, \[x\sim \sqrt{\frac{-\lambda}{2\tau}},\quad z\sim \frac{-\alpha}{\tau},\quad v\sim \sqrt{\frac{1}{-2\lambda \tau}},\quad s\sim \sqrt{\frac{2 \tau}{-\lambda}},\quad y\cdot \tau\to 0\] for some \(\alpha>0\) as \(\tau\to \infty\).
\begin{proof}
Since any trajectory converging to \(0\) that is not in the stable manifold converges exponentially to a trajectory in the centre manifold, it suffices to compute the asymptotics there. We start by computing the asymptotics for \(v=\xi_2\). Since \[\xi_2'= -\xi_2^3(-\lambda+r_2(\xi)),\] dividing through by \(\xi_2^3\) we have by L'Hopitals rule that \[v=\xi_2\sim \frac{1}{\sqrt{-2\lambda \tau}}.\] Then, Lemma~\ref{Lemma:AsymptoticsofScalarCurv} gives \(z\sim -\alpha \tau^{-1}\) for some \(\alpha>0\). The asymptotics for \(x\) follow from \(x+\lambda v=\xi_2 g(\xi)\) by dividing through by \(v=\xi_2\). The asymptotics for our original variable \(s\) are clear from \(\tau'=v^{-1} \). Finally, dividing \eqref{Eq:centremanifoldPDE} by \(\tau\) and sending \(\tau\to \infty\) we have \[0 =dh_\xi[F(\xi,h(\xi)]\tau+ \tau h(\xi)+ G(\xi,h(\xi))\tau\to 0 +\lim_{\tau\to\infty}\tau h(\xi)+0,\] where we have used the asymptotics for \(z=\xi_1,v=\xi_2\). Since \(y=\eta_2=h(\xi)_2\), this shows \(\tau \cdot y \to 0\).
\end{proof}
\end{prop}

\section{Some facts about proper group actions}\label{Appendix:ProperGroupActions}

In this appendix we collect some facts about proper group actions that are relevant to the paper. These facts are well known, however, we include proofs for convenience since they are integral to our results.

The first fact we include here is a short proof that for a proper, isometric action of \(\sans{G}\) on \((M,g)\), completeness of the orbit space \(M/\sans{G}\) is enough to conclude completeness of \(M\). Recall that the metric on \(M/\sans{G}\) is given by \[d(\sans{G}\cdot p,\sans{G}\cdot q)=\inf\{d(g\cdot p, g'\cdot q): g,g'\in \sans{G}\}.\]
Notice that the quotient map \(\pi:M\to M/\sans{G}\) is \(1\)-Lipschitz.
\begin{lemma}\label{Lemma:Completeness}
Assume that \(\sans{G}\) acts properly and isometrically on \((M,g)\) such that \(M/\sans{G}\) is a complete metric space.  Then, \((M,g)\) is complete.
\begin{proof}
Given a Cauchy sequence \(\{p_i\}\subset M\), the sequence \(\{\pi(p_i)\}\subset M/\sans{G}\) is also Cauchy. Hence, it converges to some \(p^*\in M/\sans{G}\) by completeness. Fix \(p\in \pi^{-1}(p^*)\).

Then, \(d(\sans{G}\cdot p_i,\sans{G}\cdot p  ) =d(\pi(p_i),p^*)\to 0\). For \(i\in\N\), we can find a sequence \(\{q_{ik}\}_{k=1}^\infty\subset \sans{G}\cdot p\) such that \(d(p_i,q_{ik})<d(p_i,\sans{G}\cdot p)+k^{-1}\).  The diagonal, \(q_i=q_{ii}\), then gives a Cauchy sequence in \(\sans{G}\cdot p\). Since \(\sans{G}\cdot p\) is a homogeneous space (with the induced metric), it is complete.  Hence \(q_i\to q\) for some \(q\in \sans{G}\cdot p\).
Finally, \[d(p_i,q)\le d(p_i,q_i) +d(q_i,q)\le d(p_i,\sans{G}\cdot p)+i^{-1}+d(q_i,q),\] implies \(p_i\to q\).  Hence, \(M\) is complete.

\end{proof}
\end{lemma}
Our second fact ensures that isometries between cohomogeneity one manifolds that are not homogeneous preserves the normal geodesics.

\begin{prop}\label{Prop:OrbitPreservation}
Let \(\sans{G}\) act isometrically on \((M,g)\) and \((M',g')\) with cohomogeneity one, and assume that \((M,g)\) and \((M',g')\) are not homogeneous. If \(\phi:M\to M'\) is an isometry and \(\alpha:\R\to M\) is a geodesic in \(M\) meeting the orbits orthogonally, then \(\phi\circ \alpha\) is a geodesic in \(M'\) meeting the orbits orthogonally.
\begin{proof} 
Since \(\phi\) is an isometry, \(d\phi\) maps Killing fields to Killing fields. \(\alpha'(s)\) is orthogonal to any Killing field on \((M,g)\) for all \(s\in \R\), so \((\phi\circ \alpha)'(s)=d\phi_{\alpha(s)}(\alpha'(s))\) is orthogonal to any Killing field on \((M',g')\).
\end{proof}
\end{prop}

\bibliography{Inhomogeneous_deformations_of_Einstein_Solvmanifolds}

\end{document}